\documentclass[12pt]{amsart}\topmargin=0.1in\textwidth5.9in\textheight7.85in\oddsidemargin=0.3in\evensidemargin=0.3in

\usepackage{amscd,amssymb}

\theoremstyle{plain}
\newtheorem{thm}[subsubsection]{Theorem}
\newtheorem{lem}[subsubsection]{Lemma}
\newtheorem{prop}[subsubsection]{Proposition}
\newtheorem{cor}[subsubsection]{Corollary}

\theoremstyle{definition}
\newtheorem{rk}[subsubsection]{Remark}

\newtheorem{ex}[subsubsection]{Example}

\numberwithin{equation}{section}
\newcommand{\1}{{\mathbf 1}}
\newcommand{\OO}{{\mathcal O}}

\newcommand{\h}{{\mathcal H}}

\newcommand{\LL}{{\mathcal L}}

\newcommand{\Z}{\mathbb{Z}}
\newcommand{\Q}{\mathbb{Q}}

\newcommand{\C}{\mathbb{C}}

\newcommand{\D}{\mathcal{D}}

\newcommand{\PP}{\mathbb{P}}
\newcommand{\HH}{\mathbb{H}}

\DeclareMathOperator{\Hom}{Hom}
\newcommand{\red}{{\rm red}}

\DeclareMathOperator{\id}{id}

\DeclareMathOperator{\codim}{codim}

\DeclareMathOperator{\Gr}{Gr}
\DeclareMathOperator{\Sk}{Sk}
\DeclareMathOperator{\Sing}{Sing}
\DeclareMathOperator{\Pf}{Pf}
\DeclareMathOperator{\MHM}{MHM}
\DeclareMathOperator{\MHS}{MHS}
\DeclareMathOperator{\Perv}{Perv}
\DeclareMathOperator{\rat}{rat}
\DeclareMathOperator{\Sym}{Sym}
\DeclareMathOperator{\Spec}{Spec}
\DeclareMathOperator{\vol}{vol}
\DeclareMathOperator{\Tr}{Tr}

\newcommand{\inj}{\hookrightarrow}
\newcommand{\surj}{\twoheadrightarrow}
\newcommand{\ttto}{\longrightarrow}

\begin{document}

\title[Pfaffian and the Hilbert scheme]{The Milnor fibre of the Pfaffian and the Hilbert scheme of four points on~$\C^3$}

\author[Alexandru Dimca]{Alexandru Dimca}
\address{Laboratoire J.A. Dieudonn\'e, UMR du CNRS 6621, Universit\'e de Nice Sophia Antipolis, Parc Valrose, 06108 Nice Cedex 02, FRANCE.}
\email {dimca@unice.fr}

\author[Bal\'azs Szendr\H oi]{Bal\'azs Szendr\H oi} 
\address{Mathematical Institute, University of Oxford, 24-29 St Giles', Oxford, OX1~3LB, United Kingdom.}
\email{szendroi@maths.ox.ac.uk}

\thanks{The first author is partially supported by ANR-08-BLAN-0317-02 (SEDIGA). 
The second author is partially supported by OTKA grant K61116.}

\subjclass[2000]{Primary 14C30, 14F25, 32S40; Secondary 32S55, 32S60.}

\keywords{Hilbert scheme, Pfaffian, Milnor fibration, vanishing cycles, 
intersection cohomology, mixed Hodge modules.}

\begin{abstract}
We study a natural Hodge module on the Hilbert scheme of four points on affine three-space, which categorifies the Donaldson--Thomas invariant of the Hilbert scheme. We determine the weight filtration on the Hodge module explicitly in terms of intersection cohomology complexes, and compute the E-polynomial of its cohomology. The computations make essential use of a description of the singularity of the Hilbert scheme as the degeneracy locus of the Pfaffian function.
\end{abstract}

\maketitle

\section*{Introduction}

Let $(\C^3)^{[m]}$ be the Hilbert scheme of~$m$ points on affine three-space. By 
general results, $(\C^3)^{[m]}$ is known to be a separated quasi-projective scheme, 
equipped with a proper morphism, the Hilbert--Chow morphism, to the $m$-th symmetric product
\[ \pi_m\colon (\C^3)^{[m]}\to S^m\C^3.\]
As it is well known~\cite{F}, the Hilbert scheme of points of a smooth surface is 
itself smooth and irreducible, and $\pi_m$ is a resolution of singularities 
of the symmetric product. These statements fail for threefolds; $(\C^3)^{[m]}$ is already 
singular for $m=4$, and for large $m$, it is a highly reducible and non-reduced 
scheme~\cite{BI}.

On the other hand, recent work in supersymmetric gauge theory led to a de\-scrip\-tion of the
Hilbert scheme as a degeneracy locus. A choice of affine Calabi--Yau structure on $\C^3$ induces 
an embedding of $(\C^3)^{[m]}$ into a smooth quasi-projective 
variety $M_m$ of dimension $2m^2+m$, 
which in turn is equipped with a regular function $f_m\colon M_m\to\C$, 
such that
\[ (\C^3)^{[m]} = \{ df_m=0\}\subset M_m\]
is the scheme-theoretic degeneracy locus of the function $f_m$ on the smooth variety $M_m$; 
see Proposition~\ref{prop:sup}.
The reduced space $(\C^3)^{[m]}_\red$ therefore acquires a mixed Hodge module~\cite{MHM}
\[ \Phi_m = \varphi_{f_m}\!\!\left(\Q^H_{M_m}[2m^2+m]\right)(m^2-m) \in \MHM((\C^3)^{[m]}_{\red})
\]
with underlying perverse sheaf the perverse sheaf of vanishing cycles of the function~$f_m$; 
see Section~1 for the notation, and Section~\ref{sec_phim} for an explanation for the 
shift and twist appearing here. 

The pointwise Euler characteristic of the hypercohomology of~$\Phi_m$ on 
$(\C^3)^{[m]}_{\red}$ is a~$\Z$-valued constructible function studied 
earlier in~\cite{B}. This constructible function can be integrated 
over $(\C^3)^{[m]}_{\red}$ to obtain an enumerative 
invariant, the Donaldson--Thomas invariant~\cite{T, MNOP} 
of $(\C^3)^{[m]}_{\red}$. Torus localization~\cite{MNOP, BF} shows that
the value of this invariant is the (signed) count of 
3-dimensional partitions of weight~$m$, the latter being the torus-fixed points
of the Hilbert scheme under the natural torus action. 

The module $\Phi_m$ categorifies these enumerative constructions. 
The purpose of this paper is to study the first non-trivial case $m=4$ in explicit
topological terms. We give an alternative description of $\Phi_4$, 
using the (previously known) description of the singularities of~$(\C^3)^{[4]}$ in terms of the 
Grassmannian $\Gr(2,6)$ and the Pfaffian function on skew-symmetric matrices.
We use this description to compute in Theorem~\ref{thm_E4} the compactly supported E-polynomial
\[ E_c ^{[4]}(x,y) = \sum_{p,q} \sum_k (-1)^k h^{p,q} (\HH_c^k((\C^3)^{[4]}, \Phi_4))\,x^py^q, \]
encoding the dimensions for various $p,q$ of the $(p,q)$-part of the mixed Hodge structure 
on the compactly supported hypercohomology $\HH_c^k((\C^3)^{[4]}, \Phi_4)$. The computation
makes essential use of the standard motivic properties of the 
E-polynomial, recalled in Proposition~\ref{prop_mot}.

The properties of $\Phi_m$ for general~$m$, in particular the E-polynomial $E_c^{[m]}(x,y)$, 
will be studied elsewhere~\cite{bbs} using different methods. While the general case
will present new features, the special case studied here already gives some interesting 
insight into the structure of these natural Hodge modules; compare Remark~\ref{rmk_concl}. 

In Section~\ref{sec_mhm}, we recall general facts about mixed Hodge modules, 
their E-poly\-no\-mials, and vanishing cycles. In Section~\ref{sec_pfaff}, we investigate the 
topology of the Pfaffian singularity. In Section~\ref{sec_hilb}, we 
introduce and study the 
mixed Hodge module~$\Phi_4$ on the Hilbert scheme of four points on~$\C^3$. 

\section{Mixed Hodge modules, E-polynomials and vanishing cycles}
\label{sec_mhm}
\subsection{Perverse sheaves and mixed Hodge modules}

For a reduced and separable complex algebraic variety $V$, let $\Perv(V,\Q)$ and $\MHM(V)$ denote 
respectively the abelian categories of perverse sheaves (for the middle perversity)~\cite{bbd}
and mixed Hodge modules~\cite{MHM} on~$V$. The realization functor
\[\rat\colon \MHM(V) \to \Perv(V,\Q),\]
mapping a mixed Hodge module to the underlying perverse sheaf, extends to a derived functor 
\[\rat\colon \D^b\MHM(V)\to \D^b_c(V,\Q)\]
to the derived category of constructible $\Q$-sheaves on~$V$. These triangulated categories
come equipped with the usual functors $f_*, f^*, f_!, f^!, \otimes, D$, compatibly under the
(derived) realization functor~$\rat$. 

The category of mixed Hodge modules on a point is equivalent to Deligne's category 
$\MHS$ of (polarizable) 
mixed Hodge structures, containing the canonical object $\Q^H$ of weight $(0,0)$. 
For a variety $V$, let $p\colon V\to .$ be the projection to the point, then we have
$\Q_V^H=p^*\Q^H\in\D^b\MHM(X)$; if $V$ is smooth, then in fact $\Q_V^H[\dim V]\in\MHM(V)$ 
is pure of weight $\dim V$. If $\Phi\in\D^b\MHM(V)$ is a complex of mixed Hodge modules, then
$p_*\Phi, p_!\Phi\in\MHS$, in other words the cohomologies
\[ H^*(p_*\Phi) = \HH^*(V,\Phi)\]
and
\[ H^*(p_!\Phi) = \HH_c^*(V,\Phi)\]
carry mixed Hodge structures.

\subsection{E-polynomials}
Given a variety $V$ with a bounded complex $\Phi\in\D^b\MHM(V)$ of mixed Hodge modules, let
$$E(V,\Phi;x,y) = \sum_{k,p,q} (-1)^k h^{p,q}( \HH^k(V,\Phi))x^py^q$$
be its E-polynomial and 
\[ E_c(V,\Phi;x,y) = \sum_{k,p,q} (-1)^k h^{p,q}( \HH_c^k(V,\Phi))x^py^q
\]
its compactly supported E-polynomial. For constant coefficients $\Phi=\Q_V^H$, we 
sometimes write $E(V,\Phi;x,y)=E(V;x,y)$ and $E_c(V,\Phi;x,y)=E_c(V;x,y)$.

Setting $x=y=1$ in $E, E_c$ recovers 
the Euler characteristic $\chi(V,\Phi)=\chi_c(V,\Phi)$
of (compactly supported) cohomology of $\Phi$ on $V$; see \cite[Corollary 4.1.23]{D3}
for this equality of Euler characteristics for any constructible sheaf complex.

If $\Phi[k]$ and $\Phi(k)$ denote the shift and the Tate twist of $\Phi$ by an integer $k$, 
respectively, then 
\[E(V,\Phi[k];x,y)=(-1)^{-k}E(V,\Phi;x,y)\] 
and \[E(V,\Phi(k);x,y)=(xy)^{-k}E(V,\Phi;x,y);\]
the same formulae also hold for $E_c$. 
On the other hand, Verdier duality implies that if $V$ is irreducible, then
\begin{equation} \label{VD1} 
E(V,\Phi;x,y)=E_c(V,D\Phi;x^{-1},y^{-1}).
\end{equation}
In particular, if $D(\Phi)=\Phi(n)$, then
\begin{equation} \label{VD2} 
E(V,\Phi;x,y)=x^ny^nE_c(V,D\Phi;x^{-1},y^{-1}).
\end{equation}
By~\cite[1.4]{DuS}, this applies for example to the intersection cohomology mixed 
Hodge module $\Phi=IC_V^H$, with $n=\dim V$.

\begin{prop}
\label{prop_mot} Let $V$ be a reduced variety and 
$i\colon Z\to V$ a closed inclusion of a reduced subvariety with complement 
$j\colon U\to V$. Given a complex of mixed Hodge modules $\Phi\in\D^b\MHM(V)$ on~$V$, 
we have
\[ E_c(V, \Phi;x,y) = E_c(U, j^*\Phi; x,y) + E_c(Z, i^*\Phi; x,y).
\]
\end{prop}
\begin{proof} By \cite[2.20]{MHM}, there is an exact triangle 
\[ Rj_!j^* \Phi\longrightarrow \Phi \longrightarrow i_*Li^* \Phi \stackrel{[1]}\longrightarrow
\]
in $\D^b\MHM(V)$. Taking compactly supported cohomology (compact pushforward to the point) 
gives the result. 
\end{proof}

As a consequence,  we get
\begin{cor} Let $f\colon V\to T$ be a Zariski locally trivial fibration with fibre $F$. Then 
for constant coefficients,
\[ E_c(V;x,y) = E_c(T; x,y) \cdot E_c(F; x,y).
\]
\label{cor_prod}
\end{cor}
\begin{proof} This is standard; for a product, the result simply expresses the compatibility of
the K\"unneth decomposition on compactly supported cohomology with Hodge structures. 
For a fibre bundle, stratify $T$ by strata over which $f$ is a product and use the 
previous Proposition.
\end{proof}

\subsection{Vanishing cycles} 

Let $f\colon V\to \C$ be a function on an $n$-dimensional smooth variety. Let
${}^p\varphi_f(\Q_V[n])\in\Perv(V,\Q)$ be the perverse vanishing cycle sheaf~\cite{bbd} of~$f$. 
By Saito's theory, this perverse sheaf underlies a canonical mixed Hodge module 
$\varphi_f(\Q_V^H[n])\in\MHM(V)$. The support of this module
is contained in the degeneracy locus \[ Z=\{df=0\}\subset V\]
of~$f$, so by~\cite[(2.17.5)]{MHM}, we can view $\varphi_f(\Q^H[n])\in\MHM(Z)$. 

\begin{lem} Let $V=\C^n$ be affine space, $f\colon V\to \C$ a regular function on~$V$. 
Let $W=\C^l$ be another affine space, and $g\colon V\times W\to\C$ also regular, such that
$g|_{V\times\{0\}}=0$. 
Define \[h\colon V\times W\times W^*\to\C\]
by 
\[ h(v,w,\alpha)=f(v)+g(v,w)+ \langle \alpha,w\rangle,
\]
where angle brackets denote the natural pairing between elements of $W$ and $W^*$. 
Then the projection $p: V\times W\times W^*\to V$
maps the degeneracy locus
\[ Z_h=\{ dh=0 \} \subset V\times W\times W^*
\]
of $h$ isomorphically to the degeneracy locus
\[ Z_f=\{ df=0 \} \subset V
\]
of $f$. Under this isomorphism, the vanishing cycles are related by 
\[\varphi_h\Q^H_{V\times W\times W^*}[n+2l]\cong p^*\varphi_f\Q^H_{V}[n](-l).\]
\label{lem_reduce_fn}\end{lem} 
\begin{proof} In coordinates $x=(x_1,...x_n)$ on $V$, $y=(y_1,...,y_l)$ on $W$ and dual coordinates $z=(z_1,...,z_l)$ on $W^*$, 
we have
\[ h(x, y, z) =  f(x) + g(x,y) + \sum_{j=1}^l y_jz_j\]
on $\C^{n+2l}$. 
The condition $g|_{V\times\{0\}}=0$ implies that we may write
$$g(x,y)= \sum_{j=1}^lg_j(x,y)y_j$$
for some polynomials~$g_j$.
By making the coordinate change $x'=x, ~~y'=y$ and $z_k'=z_k+g_k(x,y)$ for $k=1,\ldots, l$, 
we may assume that $g=0$.

The fact that $p$ is an isomorphism on $Z_h$ is clear, since $Z_h$ is contained in the subspace $V\times\{0\}$ given by the equations $y'=z'=0$.
 
The statement on vanishing cycles is also standard using a Thom-Sebastiani type result, see 
\cite[Section 8]{SS} for the isolated singularity case and 
\cite[Section 4]{STS} for the general case. In fact, for any $m$, the $m$-th reduced 
cohomology of the local Milnor fiber $F(h)_{(z,0,0)}$ of the $h$ at a point $(z,0,0) \in Z_h$ 
is given (as a mixed Hodge structure) by the tensor product
$$\tilde H^m(F(h)_{(z,0,0)},\Q)= \tilde H^{m-2l}(F(f)_z,\Q) \otimes  H^{2l-1}(F(A_1),\Q)$$
where $F(A_1)$ is the Milnor fibre of an $A_1$-singularity in~$\C^{2l}$. It is well known 
that $H^{2l-1}(F(A_1),\Q)=\Q(-l)$, see for instance~\cite[C26, p.243]{D2}.
\end{proof}

\begin{ex} Suppose that $V=\C^n$ and $f=\sum_{i=1}^{2k} z_i^2$ for some $k\leq \frac{n}{2}$. The
degeneracy locus is $Z=\{z_1=\ldots=z_{2k}=0\}\cong\C^{n-2k}$ carrying the 
vanishing cycle sheaf \[\varphi_f(\Q_V^H[n])\cong \Q^H_Z[n-2k](-k).\] 
The shift is given by the fact that the constructible complex underlying the vanishing cycle sheaf 
is perverse. The Tate twist comes from the fact that 
the fibre of the MHM $\varphi_{f}(\Q_V^H[n])$ at 
$z\in Z$ is given by the cohomology group $H^{2k-1}(F_z,\Q)$, where $F_z$ is the Milnor fibre 
of an $A_1$-singularity in~$\C^{2k}$. In particular, the Tate twist
$\varphi_{f}(\Q_V^H[n])(k)$ is just the constant perverse sheaf~$\Q^H_Z[\dim Z]$. 
\label{ex_A1}
\end{ex}

\section{The topology of the Pfaffian singularity} 
\label{sec_pfaff}

\subsection{The space of skew-symmetric matrices} 
Let $\Sk(2n,\C)$ be the vector space of complex skew-sym\-met\-ric $2n \times 2n$ 
matrices. It is clear that
\begin{equation} 
\label{eq1}
\dim \Sk(2n,\C)=1+2+...+(2n-1)=n(2n-1).
\end{equation}
The general linear group $G=G\ell(2n,\C)$ acts on $\Sk(2n,\C)$ by 
$g \cdot A=gAg^t$, giving rise to finitely many orbits
\begin{equation} 
\label{eq2}
V_{2k}=G A_{2k},
\end{equation}
where $k=0,1,...,n$ and 
\begin{equation} 
\label{eq3}
 A_{2k}=\left(
  \begin{array}{ccccccc}
     0_k & I_k& 0\\
     -I_k& 0_k&0 \\
      0 & 0 & 0 \\
  \end{array}
\right)
\end{equation}
is the standard skew-symmetric $2n \times 2n$ matrice of rank $2k$.
Using the techniques explained in \cite[Chapter 5]{D1}, it is easy to 
prove the following.

\begin{prop} \label{p1} There are decompositions
$$\overline  V_{2k}=\bigcup_{j=0}^{k}V_{2j},$$
with
$$\dim \overline  V_{2k}=k(4n-2k-1).$$
\label{dimprop}
\end{prop}

\begin{ex} \label{ex1} 
The non-degenerate skew matrices, i.e.~those of maximal rank $2n$, form 
a Zariski open subset, in accordance with 
\[\dim \overline  V_{2n}=n(4n-2n-1)=n(2n-1)=\dim \Sk(2n,\C).\] 
The next orbit is $V_{2n-2}$, with \[\dim \overline  V_{2n-2}=(n-1)(2n+1)=\dim \Sk(2n,\C)-1,\]
while the next orbit has dimension
$\dim \overline  V_{2n-4}=(n-2)(2n+3)=\dim \Sk(2n,\C)-6$.
\end{ex}

\subsection{The Pfaffian}\label{subsec_pf}
 It is known, see for instance \cite{L}, that there is a homogeneous 
polynomial $\Pf$ of degree~$n$, called the Pfaffian, in the entries $a_{ij}$ 
of a general matrix $A \in \Sk(2n,\C)$, such that 
\[ 
\det A= \Pf (A)^2.
\]
The hypersurface $\overline V_{2n-2}\subset\Sk(2n,\C)$ in the stratification 
of Proposition~\ref{dimprop} is the Pfaffian hypersurface $Z=\{\Pf=0\}$. 

As for any homogeneous polynomial, see \cite[Chapter 3]{D2}, there is an associated 
global Milnor fibration 
\begin{equation}\Pf: \Sk(2n,\C) \setminus Z \to \C^*,\label{fibr}\end{equation}
with central fibre the Pfaffian hypersurface~$Z$.

\begin{ex} \label{ex01}
For $n=2$, denote a general matrix $A \in \Sk(4,\C)$ by
\begin{equation} 
\label{eq02}
 A=\left(
  \begin{array}{ccccccc}
     0 & a& b & c\\
     -a& 0&d&e \\
      -b & -d & 0&f \\
      -c&-e&-f&0\\
  \end{array}
\right).
\end{equation}
Then $\Pf(A)=af-be+cd$ has an isolated singularity at the origin 
of $\Sk(4,\C)=\C^6$ of type $A_1$.
\end{ex}

\subsection{Monodromy and cohomology} \label{s02}

\begin{prop} \label{t01}
Let $F=\Pf^{-1}(1)$ be the global Milnor fibre of the 
Milnor fibration~(\ref{fibr}).

\noindent (i) The Milnor fibration~(\ref{fibr}) is algebraically trivial, thus
$$\Sk(2n,\C) \setminus Z=F \times \C^*.$$

\noindent (ii) The Milnor fibre $F$ is $4$-connected. Its Betti polynomial
$$B(F,t)=\sum_kb_k(F)t^k$$
is given by 
$$B(F,t)=(1+t^5)(1+t^9) \cdots (1+t^{4n-3}).$$
In particular, if we set $d=\dim F=2n^2-n-1$, then $b_d(F)=1$, i.e. the top 
possible non-zero Betti number of the affine smooth hypersurface $F$ is equal 
to~1.

\noindent (iii) The E-polynomial with constant coefficients of the Milnor fibre~$F$ is
$$E(F;x,y)=(1-x^3y^3)(1-x^5y^5)\cdots (1-x^{2n-1}y ^{2n-1} ).$$
\end{prop}

\begin{proof} We start by proving the easy statement (i). 
As shown for instance in~\cite{L}, under the $G$-action considered in the 
previous section, the Pfaffian satisfies
\begin{equation} 
\label{eq4}
\Pf(g\cdot A)= \det (g) \Pf(A)
\end{equation}
for any $g \in G$ and $A \in \Sk(2n,\C)$. 
It follows that the mapping \[h: F \times \C^* \to \Sk(2n,\C) \setminus Z\] 
given by 
$$ (A,t) \mapsto g(t) \cdot A,$$
with $g(t)$ the diagonal matrix $(t,1,...,1)$, 
is an isomorphism compatible with the maps to $\C^*$, given by the second 
projection and $h(A,t)\mapsto \Pf(h(A,t))=t$.

\medskip

Now we pass to statement (ii). 
The connectivity result follows from the fact that
$\codim \Sing(\Pf)=6$, see  \cite[p. 76]{D2};  
the codimension is computed using Example~\ref{ex1},
since $\Sing(\Pf)=\Sing(Z)= \overline  V_{2n-4}$. 

To compute the Betti polynomial,
let $M=\Sk(2n,\C) \setminus Z$ and note that, via the $G$-action considered 
before, $M$ is the homogeneous space $G\ell(2n,\C)/Sp(2n,\C)$. It follows 
that $M$ has the homotopy type of $M_1=U(2n)/Q(n)$, 
with $Q(n)$ the group of $n\times n$ quaternionic
orthogonal matrices. By \cite[Table I, p.493]{GHV}, 
the Betti polynomial of $M_1$ is
$$B(M_1,t)= \sum_kb_k(M_1)t^k=(1+t)(1+t^5)(1+t^9)...(1+t^{4n-3}).$$
Using the triviality of the Milnor fibration, we get 
$$B(F,t)=B(M_1,t)/(1+t)=(1+t^5)(1+t^9)...(1+t^{4n-3}).$$
This implies in particular that $b_5(F)=1$, hence the 
connectivity claim is the best possible one.

\medskip

Finally we address statement (iii) in Proposition~\ref{t01}. Using 
\cite[(6.6)]{DL} and \cite[Theorem (9.1.5)]{Del},  we have
$$E(G\ell(2n,\C);x,y)=(1-xy)(1-x^2y^2)\cdots (1-x^{2n}y^{2n}  ),$$
as well as 
$$E(Sp(2n,\C);x,y)=(1-x^2y^2 )(1-x^4y^4)\cdots (1-x^{2n}y^{2n}),$$
and hence, exactly as in \cite[(6.7)]{DL}, we have
$$E(M;x,y)=E(G\ell(2n,\C)/Sp(2n,\C);x,y)=(1-xy)(1-x^3y^3)\cdots (1-x^{2n-1}  y^{2n-1}).$$
Using the triviality of the Milnor fibration, we get 
$$E(F;x,y)=E(M;x,y)/E(\C^*;x,y)=(1-x^3y^3)\cdots (1-x^{2n-1}  y^{2n-1}).$$
This completes the proof. 
\end{proof}

Proposition~\ref{t01} shows that the Milnor fibration of the Pfaffian is 
very different from the Milnor fibration of a generic homogeneous polynomial 
of degree~$n$, where the only interesting Betti number is the top one 
and the monodromy is non-trivial.

\begin{ex} \label{ex2} The above formulas in the case $n=3$ yield the following
information on the cohomology of the Milnor fibre, where, as before, 
$(k)$ denotes Tate twist:
$$H^0(F,\Q)=\Q(0),~~H^5(F,\Q)=\Q(-3),~~H^9(F,\Q)=\Q(-5),~~H^{14}(F,\Q)=\Q(-8),$$
while the other cohomology groups of $F$ are trivial.
\end{ex}

\subsection{The perverse sheaf of vanishing cycles}
\label{s03}
We restrict in this section to the case $n=3$. Denote by 
\[X=\{d\Pf=0\}\subset Z= \{\Pf=0\}\subset \Sk(6,\C)= \C^{15}\] 
the singular locus of the $14$-dimensional Pfaffian hypersurface $Z$. 
It is easy to check that the equations $\{d\Pf=0\}$ define a reduced affine 
subvariety of $\C^{15}$ of dimension~$9$. To determine 
this variety explicitly,  
identify $A \in X= \overline  V_{2}$ with a $2$-form on~$\C^6$. 
Such a form, having rank at most $2$, can be written as
$\ell_1 \wedge \ell_2$, with $\ell_i$ linear forms on~$\C^6$.
In other words, $A$ belongs to the image of the mapping 
\[\wedge ^1\C^6 \times \wedge ^1\C^6 \to \wedge ^2\C^6 = \Sk(6,\C)=\C^{15}.\] 
As is well known, see e.g.~\cite[p. 209]{GH}, this means that
$X$ is the affine cone over the 
complex Grassmannian $Y=\Gr(2,6)$ in its Pl\"ucker embedding into $\PP^{14}$.
In particular, $X$ has an isolated singularity at the vertex $0$ of the cone. 

\begin{lem} \label{l1}
The variety $U=X \setminus \{0\}$ is simply-connected.
\end{lem}
\begin{proof}
The long exact homotopy sequence of the induced Hopf fibration on $Y$
\begin{equation} 
\label{HF}
\C^* \to U \to Y
\end{equation}
and the fact that $Y$ is simply-connected imply that $\pi_1(U)$ is 
a cyclic group (possibly trivial). The Thom long exact sequence of mixed Hodge structures
of the fibration~\eqref{HF} reads
\begin{equation} 
\label{TS}
\cdots \to H^k(Y,\Z) \to H^k(U,\Z) \to H^{k-1} (Y,\Z)(-1) \to H^{k+1} (Y,\Z) \to  \cdots
\end{equation} 
where the last morphism $H^{k-1} (Y,\Q)(-1) \to H^{k+1} (Y,\Q)$ is given by 
the cup-product by the first Chern class of the induced Hopf bundle on $Y$.
The functoriality of this long exact sequence with respect to the 
inclusion $l:Y \to \PP^{14}$ and the fact that $l^*: H^2( \PP^{14} ,\Z) \to H^2( Y ,\Z)$ is an isomorphism by \cite{La}, imply that $H^1(U,\Z)=
H^2(U,\Z)=0$. Finally this shows that $\pi_1(U)=H_1(U,\Z)=0$.
\end{proof}

Consider the perverse vanishing cycle 
sheaf ${}^p\varphi_{\Pf}(\Q[15])\in\Perv(X,\Q)$ of the Pfaffian function. 
Consider also the intersection cohomology sheaf 
\begin{equation} 
\label{eqIC}
IC_X=j_{!*}\Q_U[9]\in\Perv(X,\Q),
\end{equation}
where $j:U \to X$ is the inclusion of the smooth part, and $j_{!*}$ is 
the intermediate extension functor~\cite{bbd}. 

\begin{prop} \label{t2} Using the above notation, one has the following.

\noindent (i) On the open set $U=X \setminus \{0\}$, the two perverse sheaves
have the same restriction, namely 
\[{}^p\varphi_{\Pf}(\Q[15])|U=IC_X|U=\Q_U[9].\]
\noindent (ii) The non-zero stalks of the cohomology groups of 
${}^p\varphi_{\Pf}(\Q[15])$ at the origin are  
\[\h^k({}^p\varphi_{\Pf}(\Q[15]))_0=\Q\mbox{ for } k=0, -5, -9.\]
\noindent (iii) The non-zero stalks of the cohomology groups of 
$IC_X$ at the origin are \[\h^k(IC_X)_0=\Q \mbox{ for } k=-1, -5, -9.\]
\end{prop}
\begin{proof}
(i) The fact that $IC_X|U=\Q_U[9]$ is clear from~\eqref{eqIC}. 
To prove the other equality, let $P\in U$ be a smooth point of $X$. 
Since the critical locus $X=\{d\Pf=0\}$ is reduced and smooth near $P$,
it follows from the Morse lemma with parameters that
in appropriate analytic coordinates around $P$, the Jacobian ideal of
$\Pf$ is given by $J_{\Pf} = (x_1, \ldots, x_c)$, where $c=\codim(X)=6$,
and in these coordinates, locally $\Pf=\sum_{i=1}^c x_i^2$. 
Since the vanishing cycle sheaf of this Morse
function is of rank one, concentrated in the appropriate degree,
it follows that there exists a rank one local system $\LL$ on $U$ 
such that ${}^p\varphi_{\Pf}(\Q[15])=\LL[9].$
In view of Lemma \ref{l1}, necessarily $\LL=\Q_U$, completing
the proof of (i).

Statement (ii) follows directly from Proposition~\ref{t01}~(ii) above, upon noting
\[\h^k({}^p\varphi_{\Pf}(\Q[15]))_0=\tilde  H^{14+k}(F,\Q).\]

To prove (iii), 
recall that since $0$ is an isolated singularity of $X$, 
$$IC_X=\tau_{\leq -1}(Rj_*\Q_U[9]),$$
see \cite[Prop 5.2.10, Proof of Prop. 5.4.4]{D3}. 
This implies that 
\[\h^k(IC_X)_0=\h^k(Rj_*\Q_U[9] )_0\] for $k\leq -1$ and 
$\h^0(IC_X)_0=0$.

On the other hand, $\h^k(Rj_*\Q_U[9] )_0=H^{k+9}(U,\Q)$ since $U$ has the 
same homotopy type as the link of the singularity $(X,0)$. In the long 
exact sequence~\eqref{TS} associated to the fibration~\eqref{HF}, 
the last morphism $H^{k-1} (Y,\Q)\to H^{k+1} (Y,\Q)$ is 
injective for $k \leq \dim Y=8$ 
by the Hard Lefschetz Theorem, since it is given by the cup-product by the 
first Chern class of the induced Hopf bundle on $Y$, a non-zero multiple of 
the K\"ahler class of $Y$. This implies that, for $k\leq 9$, we have 
$$b_k(U)=b_k(Y)-b_{k-2}(Y).$$
The Betti polynomial of the Grassmannian $Y=\Gr(2,6)$ is well known:
$$B(Y,t)=\frac{(1-t^{10})(1-t^{12})}{(1-t^{2})(1-t^{4})}.$$
It follows that the only non-zero Betti numbers $b_k(U)$ for $k \leq 9$ 
are
$$b_0(U)=b_4(U)=b_8(U)=1.$$
This completes the proof of~(iii).

\end{proof}

\begin{rk} \label{MHSU} For future use, note that the exact sequence~\eqref{TS} of mixed Hodge 
structures and the fact that the cohomology groups $H^{2k}(Y,\Q)$ have pure Hodge type $(k,k)$ 
for all $k$, imply that
$H^0(U,\Q)=\Q(0),~~H^4(U,\Q)=\Q(-2),~~H^8(U,\Q)=\Q(-4),~~H^9(U,\Q)=\Q(-5),~~H^{13}(U,\Q)=\Q(-7),~~H^{17}(U,\Q)=\Q(-9).$

\end{rk}

\begin{rk} \label{conicrk}
Since $X$ is a cone and both sheaves ${}^p\varphi_{\Pf}(\Q[15])$ and $IC_X$ are constructible with 
respect to the obvious stratification $(\{0\},U)$ of $X$, it follows that for any $k$, 
$$ \HH^k(X, {}^p\varphi_{\Pf}(\Q[15]))=\HH^k(B \cap X, {}^p\varphi_{\Pf}(\Q[15]))=\h^k( {}^p\varphi_{\Pf}(\Q[15]))_0,$$
and similarly
$$ \HH^k(X, IC_X)=\HH^k(B \cap X, IC_X)=\h^k(IC_X)_0.$$
Here $B$ is a (small) ball centered at the origin, the first isomorphisms come from the conic 
structure of $X$ and the ones involving the stalk at $0$ from \cite[Corollary 4.3.11]{D3}. 
\end{rk}

Denote by $i:0\to X$ the inclusion of the singular point.

\begin{thm} In the category $\Perv(X,\Q)$ of perverse sheaves on $X$, there is a three-step
filtration on ${}^p\varphi_{\Pf}(\Q[15])$, with the quotients being $i_*\Q_0$,
$IC_X$ and $i_*\Q_0$ respectively. 
\label{t3}
\end{thm}
\begin{proof} Following the arguments of~\cite[pp. 134-135]{D3}, 
Theorem~\ref{t2} (i) above shows that there exist a commutative
diagram 
\[\begin{array}{ccccc} {}^p Rj_!\Q_U[9] & \stackrel{\alpha}\ttto & {}^p\varphi_{\Pf}(\Q[15]) & \ttto & {}^p R j_*\Q_U[9] \\
\parallel &&&& \parallel\\
{}^p Rj_!\Q_U[9] & \surj & IC_X & \inj & {}^p Rj_*\Q_U[9].\end{array}\]
Let $A$ denote the image of the map $\alpha$. Then we 
obtain the diagram 
\[\begin{array}{ccccc} {}^p Rj_!\Q_U[9] & \surj & A & \ttto & {}^p R j_*\Q_U[9] \\
\parallel &&&& \parallel\\
{}^p Rj_!\Q_U[9] & \surj & IC_X & \inj & {}^p Rj_*\Q_U[9].\end{array}\]
The definition of $A$, and an easy diagram chase in this diagram give 
two short exact sequences
\begin{equation}\label{seq1} 0\to A \to {}^p\varphi_{\Pf}(\Q[15])\to Q_1 \to 0\end{equation}
and
\begin{equation}\label{seq2} 0\to Q_2\to A \to IC_X \to 0\end{equation}
in the abelian category of perverse sheaves on~$X$. 
Again by Theorem~\ref{t2} (i), the perverse sheaves $Q_j$ are supported
on the singular point $0$, and so they are of the form $i_*V_j$ for 
$\Q$-vector spaces~$V_j$. Taking hypercohomology at~$0$, using 
$\h^{k}(i_*V_j)_0=0$ for $k\neq 0$, 
as well as the results of  
Theorem~\ref{t2} (ii)-(iii), we get two short exact
sequences of vector spaces
\[ 0\to\h^{0}(A)_0 \to \Q \to V_1 \to 0,
\]
and
\[ 0\to\h^{-1}(A)_0 \to \Q \to V_2 \to 0.
\]
Clearly $V_1=0$, $V_2=0$ are impossible, so they are both one-dimensional
$\Q$-vector spaces. Hence the sequences~\eqref{seq1}-\eqref{seq2} exhibit 
the claimed filtration on the vanishing cycle 
sheaf ${}^p\varphi_{\Pf}(\Q[15])$.
\end{proof}

Note that by \cite[Prop. 4.2.10, Prop. 3.3.7, Prop. 5.2.9]{D3} respectively, 
the sheaves ${}^p\varphi_{\Pf}(\Q[15])$, $IC_X$ and $i_*\Q_0$
are all self-dual under Verdier duality. This is reflected by the
fact that the series of composition factors in the filtration of the sheaf of 
vanishing cycles is palindromic.

\subsection{The vanishing cycle Hodge module and its cohomology} 

In this section, we lift the results of the previous section to the category~$\MHM(X)$
of mixed Hodge modules on~$X$. 
By Saito's general theory, the perverse sheaf ${}^p\varphi_{\Pf}(\Q[15])$ underlies
a canonical mixed Hodge module $\varphi_{\Pf}(\Q^H[15])\in\MHM(X)$.
 
\begin{thm} \label{main1} In the category $\MHM(X)$ of mixed Hodge modules on $X$, 
the weight filtration
on $\varphi_{\Pf}(\Q^H[15])$ has quotients $i_*\Q^H_0(-7)$ in weight $14$, 
$IC^H_X(-3)$ in weight $15$ and $i_*\Q^H_0(-8)$ in weight $16$ respectively. Here $IC^H_X$ and $\Q^H_0$ are
the natural mixed Hodge modules over the perverse sheaves $IC_X$ and $\Q_0$ on $X$ and $0$, 
respectively. 
\end{thm} 
\begin{proof} The graded quotients of the weight filtration on $\varphi_{\Pf}(\Q^H[15])$
are pure Hodge modules, which in turn decompose into a direct sum of objects with strict support. 
Hence the filtration found in Theorem~\ref{t3} is necessarily the filtration by perverse 
sheaves underlying the weight filtration. The proof of Theorem~\ref{t3} also shows that the
graded pieces are as claimed, once we replace the isomorphism
$${}^p\varphi_{\Pf}(\Q[15])|U=IC_X|U=\Q_U[9]$$
from Proposition \ref{t2} by the mixed Hodge module isomorphism
$$\varphi_{\Pf}(\Q^H[15])|U=IC_X(-3)|U=\Q^H_U(-3)[9].$$
For the Tate twist, see Example~\ref{ex_A1}.
 
Next, one may use the exact sequence \eqref{seq2} and Remark \ref{MHSU} to identify
$$Q_2=\h^{-1}(IC_X)_0(-3)=H^8(U,\Q)(-3)=\Q(-7).$$
On the other hand, the exact sequence \eqref{seq1} gives
$$Q_1=\h^0({}^p\varphi_{\Pf}(\Q^H[15]))=H^{14}(F,\Q)=\Q(-8)$$
according to Example \ref{ex2}.

Finally, it is well known that for any irreducible variety $V$, the intersection complex
$IC_V^H$ is pure of weight $n=\dim V$, see for instance  \cite[Lemma 14.15]{PS}.
\end{proof}

\begin{thm} The E-polynomials of the mixed Hodge Module 
$\varphi_{\Pf}(\Q^H[15]))$ are given by 
\[ E(X,\varphi_{\Pf}(\Q^H[15]); x,y)= (xy)^3((xy)^5-(xy)^2-1)
\]
and
\[ E_c(X,\varphi_{\Pf}(\Q^H[15]); x,y)= (xy)^7(1-(xy)^{3}-(xy)^{5}).
\]
\label{cor_E}
\end{thm}
\begin{proof} 
Using Remark \ref{conicrk} we get the following isomorphisms of MHS
$$ \HH^k(X, \varphi_{\Pf}(\Q^H[15]))=\h^k( \varphi_{\Pf}(\Q^H[15]))_0= \tilde  H^{14+k}(F,\Q).$$
The formula for $E(X,\varphi_{\Pf}(\Q^H[15]); x,y)$ follows then from Example \ref{ex2}.
To compute the E-polynomial with compact supports, note that $\varphi_{\Pf}(\Q^H[15])$ 
contains as a direct summand $\varphi_{\Pf,1}(\Q^H[15])$, on which the semisimple part of the 
monodromy acts trivially. Since the underlying perverse sheaf ${}^p\varphi_{\Pf}(\Q[15])$ 
admits no non-trivial direct summand, we must have 
\[\varphi_{\Pf,1}(\Q^H[15])=\varphi_{\Pf}(\Q^H[15]).\] 
Now using \cite[(2.6.2)]{MHM} and the 
identification $IC_{\C^{15}}^H=\Q^H[15]$, we get
$$D(\varphi_{\Pf}(\Q^H[15]))=\varphi_{\Pf}(D(\Q^H[15]))=\varphi_{\Pf}(\Q^H[15])(15).$$
Hence~\eqref{VD2} applies.
\end{proof}

\section{The Hilbert scheme of four points on affine three-space}
\label{sec_hilb}

\subsection{The superpotential description}\label{sec:sup}

Let $T$ be the three-dimensional space of linear functions on $\C^3$, so 
that $\C^3=\Spec\Sym^\bullet T$. Fix an isomorphism $\vol: \wedge^3 T\cong \C$; this 
corresponds to choosing a holomorphic volume form (Calabi--Yau form) on $\C^3$. 
We start by recalling the de\-scrip\-tion of the Hilbert scheme as a degeneracy locus 
from \cite[Proof of Thm.~1.3.1]{Sz} and \cite[Prop.3.8]{seg}. 

\begin{prop} The pair $(T,\vol)$ defines an embedding of the Hilbert scheme 
$(\C^3)^{[m]}$ into a smooth quasi-projective variety $M_m$, which in turn is equipped 
with a regular function $f_m\colon M_m\to\C$, such that
\begin{equation} (\C^3)^{[m]} = \{ df_m=0\}\subset M_m\label{eq_emb}\end{equation}
is the scheme-theoretic degeneracy locus of the function $f_m$ on $M_m$.
\label{prop:sup}
\end{prop}

\begin{proof} A point $[Z]\in (\C^3)^{[m]}$ corresponds to an embedded
$0$-dimensional subscheme $Z\hookrightarrow\C^3$ of length $m$, in other words 
to a quotient $\OO_{\C^3}\to \OO_Z$ with $H^0(\OO_Z)$ of dimension $m$. 
Fixing an $m$-dimensional complex vector space~$W_m$, the data defining a cluster 
consists of a linear map $T\otimes W_m\to W_m$, subject to the condition that 
the induced action of the tensor algebra of $T$ factors through an action of the symmetric 
algebra $\Sym^\bullet T$, and a vector $1\in W_m$ which generates $W_m$ under the action.

Let
\[ U_m\subset \Hom(T\otimes W_m, W_m) \times W_m\]
denote the space of maps with cyclic vector, but without the symmetry condition.
As proved in~\cite[Lemma 1.2.1]{Sz}, the action of $G\ell(W_m)$ on $U_m$~is free,
and the quotient 
\[ M_m = U_m/G\ell(W_m)\]
is a smooth quasiprojective GIT quotient. 

Finally consider the map 
\[ \phi \mapsto  \Tr\left(\wedge^3\phi\right),\]
where $\wedge^3\phi:\bigwedge^3 T\times W_m\to W_m$ and
we use the isomorphism $\vol$ before taking the trace on~$W_m$.
It is clear that this map descends to a regular map $f_m\colon M_m\to\C$. 
As proved in~\cite[Prop.3.8]{seg}, the equations $\{df_m=0\}$ are just the the equations
which say that the action factors through the symmetric algebra. 
As proved by~\cite{nak} (in dimension 2, but the proof generalizes), 
the scheme cut out by these equations is precisely the moduli 
scheme representing the functor of $m$ points on $\C^3$. Thus, as a scheme, 
\[ (\C^3)^{[m]}=\{df_m=0\} \subset M_m.
\]
\end{proof}

\begin{rk} In explicit terms, fixing a basis of~$V$, $\C[x,y,z]$ acts on $W_m$ by
a triple of commuting matrices $X,Y,Z$. The map $f_m$ on triples of matrices is given by
\[ (X,Y,Z)\mapsto {\rm Tr} [X,Y]Z.\]
\end{rk}

\subsection{The Hodge module $\Phi_m$ on the Hilbert scheme}
\label{sec_phim}

As a consequence of Proposition~\ref{prop:sup}, 
the reduced space $(\C^3)^{[m]}_\red$ acquires a mixed Hodge module 
\[\Phi_m=\varphi_{f_m}\!\!\left(\Q^H[\dim M_m]\right)(m^2-m),\] with underlying perverse sheaf
the perverse sheaf of vanishing cycles of the function~$f_m$. 
By~\cite[1.2]{B}, the pointwise Euler characteristic
of $\Phi_m$ is the Behrend function of~$(\C^3)^{[m]}$, and hence 
the Euler characteristic over the whole
space computes its Donaldson--Thomas invariant~\cite{MNOP}. 
In this sense, this Hodge module categorifies the Donaldson--Thomas invariant of the
Hilbert scheme. 

The shift and Tate twist appearing in the definition of $\Phi_m$ are not essential, and 
are only inserted for cosmetic purposes. The shift by $\dim M_m$ simply means that 
the underlying constructible complex is a perverse sheaf; thus, it puts the object
in the abelian category $\MHM((\C^3)^{[m]}_\red)$ and not just in its derived category. 
The Tate twist cancels certain Tate twists arising from the embedding in $M_n$. 
To illustrate this, note that for $m\leq 3$, the Hilbert scheme $(\C^3)^{[m]}$ 
is nonsingular. It is moreover easy to see, using the argument in the proof of
Proposition~\ref{t2}(i), that the Hodge module~$\Phi_m$ is the trivial sheaf 
$\Q^H[\dim(\C^3)^{[m]}]$ in these cases; the Tate twist in the definition of~$\Phi_m$ 
cancels the Tate twist arising from Example~\ref{ex_A1}. See Remark~\ref{rmk_concl} for
further discussion of this point. 
The forthcoming~\cite{bbs} will compute the cohomology of
all the Hodge modules $\Phi_m$, making contact with the refined topological vertex
calculations of~\cite{ikv}; the Tate twist will be important there to get a closed
formula. 

Recent work of Kontsevich and Soibelman~\cite{ks} constructs a motivic generalization 
of Donaldson--Thomas theory in great generality, using local data. So their point of
view is slightly different; the Hodge module $\Phi_m$ exist globally on the moduli
space~$(\C^3)^{[m]}$. It seems to be an interesting question when does there exist a global 
mixed Hodge module on a moduli space ${\mathcal M}$ admitting a symmetric 
perfect obstruction theory~\cite{BF}, whose pointwise Euler characteristic is the Behrend
function of~${\mathcal M}$, and if so, then how canonical is it. 

\subsection{The Hodge module $\Phi_4$ on the Hilbert scheme of four points}
\label{sec_strat}
For the rest of the paper, consider the case $m=4$ of four points on affine 3-space. 
By~\cite{katz}, the scheme $(\C^3)^{[4]}$ is irreducible of dimension 12, singular 
along the locus~$S_4$ of length-four subschemes of $\C^3$ given 
by the squares of maximal ideals of points. 
Let $N_4=(\C^3)^{[4]}\setminus S_4$  be the nonsingular part.

Fix an affine structure on $\C^3$ once again, then there is a stratification 
\[  L_4\subset P_4 \subset (\C^3)^{[4]}
\]
defined as follows: $P_4$ is the locus of subschemes scheme-theoretically contained in a 
linear hyperplane in $\C^3$; $L_4$ is the locus of subschemes scheme-theoretically contained 
in a line. The complement 
\[ V_4 = (\C^3)^{[4]} \setminus P_4\]
is the open subset of all colength-four ideals $I\lhd \Sym^\bullet T$ such that 
as a vector space, $\Sym^\bullet T/I\cong \1\oplus T$, where $\1$ is the one-dimensional
space generated by the cyclic vector~$1$. This open subset $V_4$ contains the 
singular locus $S_4$ of $(\C^3)^{[4]}$. 

\begin{prop}\label{prop_rest}
\noindent (i) {\rm (compare \cite{katz, lee})}
The open subset $V_4$ of $(\C^3)^{[4]}$ is affine.
There is a product decomposition $V_4= \C^3\times X$, where
$X\subset\C^{15}$ is the cone over $\Gr(2,6)$ in its Pl\"ucker embedding in $\PP^{14}$.

\noindent (ii) The restriction of $\Phi_4$ to $V_4$ is 
given by $p_2^*(\varphi_{\Pf}(\Q^H[15]))[3](3)$, where $p_2\colon V_4\to X$ 
is the second projection.

\noindent (iii) The restriction of $\Phi_4$ to the nonsingular part 
$N_4=(\C^3)^{[4]}\setminus S_4$ 
is the rank-one sheaf $\Q_N^H[12]$. 
\end{prop}

\begin{proof} First we prove (i) and (ii) together. 
As remarked above, on the open set $V_4\subset(\C^3)^{[4]}$,
there is a decomposition $W_4\cong \1\oplus T$. Unwinding the definitions 
of the map $\phi$ from Proposition~\ref{prop:sup}, we can then write
\[\phi=\left(\begin{array}{cc} 0 & \phi_1\\ \id_T & \phi_2\end{array}\right)\colon T\otimes \left( \1\oplus T\right)\to \1\oplus T,
\]
so the data is equivalent to a pair of maps
\[ (\phi_1,\phi_2) \colon T\otimes T\to \1\oplus T\]
leading to an embedding
\[ V_4\subset\Hom(T\otimes T, \1\oplus T).\]
Note there is no $G\ell$-action left here, since we completely 
rigidified $W_4$ using $T$. 
This is the restriction of the embedding~\eqref{eq_emb}, mapping the open 
set $V_4\subset(\C^3)^{[4]}$ into a $36$-dimensional affine space. Thus $V_4$ is 
indeed affine. 

Write $\phi_2=\phi_2^+ + \phi_2^-$, where $\phi_2^+\in\Hom(\Sym^2 T,T)$ and  
$\phi_2^-\in\Hom(\bigwedge^2 T,T)$. Then a short computation shows
\[ \Tr\left(\wedge^3\phi\right) = \Tr\left(\wedge^3\phi_2^+\right) + g(\phi_2^+, \phi_2^-)+ 2\langle \phi_1, \phi_2^-\rangle.
\]
Here $g$ is some cubic function, which we will not need explicitly,
and the angle brackets denote the natural pairing between $\Hom(T\otimes T,\1)$ and 
\[\Hom\left(\bigwedge^2 T,T\right)\cong \Hom(T^*, T) \cong T\otimes T,\]
where the first isomorphism is induced by $\vol$. 
The conditions of Lemma \ref{lem_reduce_fn} can be checked to hold; so by that Lemma, 
the vanishing cycle sheaf of $f_4$ on $V_4$ can 
be computed from the embedding \[V_4\subset\Hom(\Sym^2 T, T)\] as the degenerecy locus
of the function $\Tr\left(\wedge^3\phi_2^+\right)$. 

Next, apply Lemma~\ref{lem_linalg} below, with $S=T^*$. Using a trivialization 
$\bigwedge^3 T^*\cong \1^*\cong \1$ induced by $\vol\colon\bigwedge^3 T\cong \1$, 
we have 
\begin{eqnarray*} \Hom(\Sym^2 T, T) & = & \Spec\Sym^\bullet (\Sym^2 T^*, T^*) \\
& \cong & \Spec\Sym^\bullet\left( T\oplus\bigwedge^2\Sym^2 T\right) \\
& \cong & \C^3 \times \bigwedge^2\Sym^2 T^*.
\end{eqnarray*}
It can be checked by explicit calculation that under this isomorphism, 
the function $\Tr\left(\wedge^3\phi_2^+\right)$ only depends on the projection to the second
factor, leading to an isomorphism  
\[ V_4\cong \C^3 \times X\subset \C^3  \times \bigwedge^2\Sym^2 T^*,\]
and the function on $\bigwedge^2\Sym^2 T^*$ is the Pfaffian of Section~\ref{subsec_pf}. 
This proves (i). As for the vanishing cycle 
sheaf, we finally have
\[\Phi_4|_{V_4} \cong p_2^*(\varphi_{\Pf}(\Q^H[15]))[3](3),\] 
where $p_2:V_4=\C^3 \times X \to X$ is the second projection.
Here the additional shift $[3]$ is due to the fact that the relative dimension of $p_2$ 
is 3; the Tate twist $3=12-9$ arises from the Tate twist in the definition of~$\Phi_4$,
and the reduction from the $36$-dimensional embedding space $(\C^3)^{[4]}\subset M_4$ to the 
$18$-dimensional embedding space $\Hom(\Sym^2 T, T)$ using Lemma~\ref{lem_reduce_fn}.
This concludes the proof of~(ii). 

To prove (iii), first of all note that the nonsingular part $(\C^3)^{[4]}\setminus S_4$ has
a Zariski open subset $V_4\setminus S_4$ which is isomorphic by~(i) 
to $U\times\C^3$, where $U$ is the nonsingular part of $X$. Since $U$ is simply connected by 
Lemma \ref{l1}, so is~$N_4$. Thus, the restriction of~$\Phi_4$ to 
$N_4=(\C^3)^{[4]}\setminus S_4$ is $\Q_{N_4}^H[12]$, since the transversal singularity is $A_1$
embedded in $\C^{24}$, exactly as in the proof of Theorem~\ref{main1}; the Tate twists 
by~$12$ in the definition of~$\Phi_4$ and $-12$ computed in Example~\ref{ex_A1} cancel. 
\end{proof}

\begin{cor} In the category $\MHM((\C^3)^{[4]})$ of mixed Hodge modules on $(\C^3)^{[4]}$, 
the weight filtration on $\Phi_4$ has quotients $j_*\Q^H_{S_4}[3](-4)$ in weight $11$, 
$IC^H_{(\C^3)^{[4]}}$ in weight $12$ and $j_*\Q^H_{S_4}[3](-5)$ in weight $13$ respectively, where
$j\colon S_4\cong\C^3 \hookrightarrow (\C^3)^{[4]}$ denotes the inclusion of the singular locus.
\end{cor} 
\begin{proof} This follows from Proposition~\ref{prop_rest}, in light of Theorem~\ref{main1}, by
applying the functor $p_2^*(-)[3](3)$ to all the objects in Theorem~\ref{main1}. Note that
$$p_2^*(IC_X^H)[3]\cong IC_{V_4}^H;$$
since both mixed Hodge modules here are pure of weight 12, it is enough 
to check that they are isomorphic on the perverse sheaf level, where it follows for 
instance using \cite[Thm. 3.2.13 (ii), Thm. 3.2.17 (iii) and Prop. 5.2.10]{D3}.
\end{proof}

The proof of Proposition~\ref{prop_rest} used the following

\begin{lem} {\rm (compare \cite[Lemma 1.6]{katz})} 
For a three-dimensional vector space $S$, there is a canonical isomorphism 
\[\Hom(\Sym^2 S, S)\cong  S^*\oplus\bigwedge^2\Sym^2 S^*\otimes \bigwedge^3 S.\]
\label{lem_linalg}
\end{lem}
\begin{proof} For any vector space $S$, there are canonical inclusions
\begin{eqnarray*} 
S^* & \to & \Hom(\Sym^2 S, S)\\
\alpha &\mapsto &\left((s_1,s_2)\mapsto \alpha(s_1)s_2 + \alpha(s_2)s_1\right)
\end{eqnarray*}
and
\[\setlength{\arraycolsep}{.4\arraycolsep}
\begin{array}{rcl}
\displaystyle\bigwedge^2\Sym^2 S^*\otimes\bigwedge^3 S & \to & \Hom(\Sym^2 S, S)\\
\displaystyle(\omega_1\wedge\omega_2)\otimes(s_1\wedge s_2\wedge s_3)&\mapsto & \displaystyle\left(\!(t_1,t_2)\mapsto 
\sum_{\sigma\in S_3} (-1)^{{\rm sign}(\sigma)}\omega_1(t_1s_{\sigma(1)})\omega_2(t_2s_{\sigma(2)})s_{\sigma(3)}\!\right)
\end{array}\]
where $t_i, s_i\in S$, $\omega_i\in\Sym^2 S^*$. Writing these maps out on a basis, one
checks that the images intersect trivially. For $\dim S=3$, 
the Lemma follows by a dimension count. 
\end{proof}

\subsection{The E-polynomial of $\Phi_4$ on the Hilbert scheme of four points}

As in the Introduction, let 
\[ E_c^{[4]}(x,y) = E_c((\C^3)^{[4]},\Phi_4; x,y)
\]
be the E-polynomial of compactly supported cohomology of the Hodge module~$\Phi_4$
on~$(\C^3)^{[4]}$.

\begin{thm} We have 
\label{thm_E4}
\[ E_c^{[4]} (x,y)  = (xy)^6((xy)^{6} + (xy)^{5} + 3 (xy)^{4}+ 3 (xy)^3+ 3 (xy)^2+ (xy)+ 1).\]
\end{thm}
\begin{proof} The computation makes use of the motivic nature of the compactly supported
E-polynomial. We use the stratification  $L_4\subset P_4 \subset (\C^3)^{[4]}$ 
from Section~\ref{sec_strat}. 

The contribution from the open stratum $V_4=(\C^3)^{[4]}\setminus P_4$ is 
\begin{eqnarray*} 
E_c(V_4, \Phi_4|_{V_4}; x,y) & = &  E_c(V_4,  p_2^*(\varphi_{\Pf}(\Q[15]))[3](3))  ; x,y) \\
&=& E(V_4, D( p_2^*(\varphi_{\Pf}(\Q[15]))[3](3))  ; x^{-1},y^{-1}) \\
&=& E(V_4,  p_2^!(\varphi_{\Pf}(\Q[15])(15))[-3](-3)  ; x^{-1},y^{-1}) \\
& = & -E(V_4,  p_2^*(\varphi_{\Pf}(\Q[15])(15))  ; x^{-1},y^{-1})\\
&=&-E(X, \varphi_{\Pf}(\Q[15])(15) ; x^{-1},y^{-1})\\
&=&-E_c(X, D(\varphi_{\Pf}(\Q[15])(15))  ; x,y)\\
&=&-E_c(X,\varphi_{\Pf}(\Q[15]);x,y)\\
&=&(xy)^{7}((xy)^5 + (xy)^{3} - 1).\end{eqnarray*}
Here, apart from standard properties of pullback and duality, we used 
Proposition~\ref{prop_rest}(ii)
for the first equality, the 
identity $p_2^!=p_2^*(3)[6]$ from \cite[(3.5.1) and (3.5.2)]{MHP} for
the fourth, and Theorem~\ref{cor_E} for the last one. 

Next, look at the smallest stratum~$L_4$ of clusters lying in a line. By 
Proposition~\ref{prop_rest}(iii), the restriction of $\Phi_4$ to this stratum 
is~$\Q^H_{L_4}[12]$. Associating to a cluster in $L_4$ the line it is contained 
in defines a map $l\colon L_4\to {\rm Lines}(\C^3)$
to the space ${\rm Lines}(\C^3)$ of all lines in $\C^3$. 
The map $l$ is easily seen to be a Zariski locally trivial
fibration, with fibre the Hilbert scheme of four points on the affine line. 
The space ${\rm Lines}(\C^3)$ of all lines in $\C^3$ further fibres 
over $\PP^2$, with fibre~$\C^2$. The Hilbert 
scheme of four points on the affine line is simply $\C^4$. 
Putting all this together, using Corollary~\ref{cor_prod} repeatedly, 
we get, for constant coefficients, 
\[ E_c(L_4; x,y) = (xy)^4 (xy)^2 (1+xy+(xy)^2).\]
Thus the contribution from $L_4$ to $E_c ^{[4]} (x,y)$ is
\[ E_c(L_4, \Phi_4|_{L_4}; x,y) = (xy)^6(1+xy+(xy)^2).\]

Finally, look at the stratum $P_4\setminus L_4$ of strictly planar clusters. By 
Proposition~\ref{prop_rest}(iii) again, the restriction of $\Phi_4$ to this stratum is still 
a shift of the trivial sheaf~$\Q^H_{P_4\setminus L_4}$.
The variety $P_4\setminus L_4$ has a Zariski locally trivial
fibration over the space ${\rm Planes}(\C^3)$ of all planes in $\C^3$, with fibre the Hilbert 
scheme of four points on the affine plane which are not collinear. 
The space ${\rm Planes}(\C^3)$ of all planes in $\C^3$ further fibres 
over $\PP^2$, with fibre~$\C$. On the other hand, the E-polynomial with constant coefficients 
of the Hilbert scheme of four points on the affine plane can be computed from 
G\"ottsche's formula~\cite{goe} to be
\[E((\C^2)^{[4]}; x,y)= (xy)^5 + 2 (xy)^6 + (xy)^7 + (xy)^8.\]
We need to subtract from this the contribution from the locus of clusters on a plane contained 
in a line, which can be computed by fibering it over the space of all lines in a plane. 
Putting all this together, we obtain
\begin{eqnarray*} E_c(P_4\setminus L_4; x,y) &=&\left((xy)^5 + 2 (xy)^6 + (xy)^7 + (xy)^8 - (xy)^4 xy (1+xy)\right)\cdot \\ && \cdot xy\cdot (1+xy+(xy)^2).\end{eqnarray*}
Collecting terms, the contribution from $P_4 \setminus L_4$ to $E_c ^{[4]} (x,y)$ is
$$E_c(P_4\setminus L_4, \Phi_4|_{P_4\setminus L_4}; x,y)=(xy)^7(1+xy+(xy)^2)^2.$$
Summing the three contributions gives the answer stated in the Theorem. 
\end{proof}

\begin{rk} For $x=y=1$, we get  $E_c ^{[4]} (1,1)=13$, which is indeed the Donaldson--Thomas 
invariant of the Hilbert scheme of four points on three-space~\cite{MNOP, BF}. As discussed in the 
Introduction, the latter is the number of 3-dimensional partitions of $4$, the number
of torus-fixed points on the Hilbert scheme. 
The formula of Theorem~\ref{thm_E4} refines this enumerative count to a polynomial. 

As a result of the Tate twist in the definition of~$\Phi_m$,
the expression for $E_c^{[4]}$ looks like the E-polynomial of a smooth
variety of dimension~$12$. Being a sum of positive powers of~$xy$, 
the formula could in fact correspond to a Hodge structure which is an effective sum
of Tate motives. This indeed happens in the much simpler cases $m\leq 3$: the Hilbert scheme
$(\C^3)^{[m]}$ in these cases is smooth, the module $\Phi_m$ is simply the constant module, and
the Hilbert scheme itself is just a union of pieces each of which is affine space
of some dimension. The E-polynomials in these cases are computed in~\cite{ch}.
However, the proof shows that for $m=4$, this cannot be true, and
one cannot recover the cohomology from the E-polynomial, because of cancellation between
terms in different degrees. 

A different way to compute the E-polynomial $E_c ^{[4]} (x,y)$ would be by torus 
localization, as in the proof of G\"ottsche's formula for $\C^2$ given in \cite[5.2]{nak}, 
the Euler characteristic computation of~\cite{BF} on the Hilbert scheme of $\C^3$, and
the computation of the Poincar\'e polynomial of the resolution of singularities of $(\PP^3)^{[4]}$
in~\cite{katz}.
As an ideal scenario, one could have hoped in fact that each of the 13 torus-fixed points 
contributes a monomial of $xy$ to the E-polynomial. But this turns out to be false. 
Of the 13 torus-fixed points corresponding to 3-dimensional partitions of 4, all but one
are planar and therefore lie in the smooth locus of the Hilbert scheme. The contribution of 
these points to the E-polynomial is easily computed using formulae from~\cite{MNOP},
a computer, and an appropriate one-dimensional subtorus of $(\C^*)^3$. The result is 
that the smooth torus-fixed points contribute
\[(xy)^{12} + 2 (xy)^{11} + 3 (xy)^{10}+ (xy)^9+ 3 (xy)^8+ (xy)^7+ (xy)^6.
\]
Comparing this with the final formula, we see that the last remaining fixed
point, the square of the maximal ideal at the origin which lies in the singular locus, 
must contribute $2(xy)^9-(xy)^{11}$. 
\label{rmk_concl}\end{rk}

\section*{Acknowledgements} The authors would like to thank
Barbara Fantechi and Lothar G\"ottsche for pointing out 
Proposition~\ref{prop_rest}(i), as well as Tom Bridgeland, Patrick Brosnan, 
Stefan Papadima, Morihiko Saito and 
Duco van Straten for helpful remarks and correspondence.

\end{document}